\documentclass[12pt]{amsart}
\usepackage{amsmath,amstext,amsthm,amssymb,amsxtra}
\usepackage[colorlinks,citecolor=red,pagebackref,hypertexnames=false]{hyperref}
\usepackage[backrefs]{amsrefs}

\usepackage{pgf,tikz}
 \allowdisplaybreaks \swapnumbers

\theoremstyle{plain} 
\newtheorem{lemma}[equation]{Lemma}
\newtheorem{proposition}[equation]{Proposition}
\newtheorem{theorem}[equation]{Theorem}

\newtheorem{definition}[equation]{Definition}

\numberwithin{equation}{section}

\begin{document}

\title[Directional discrepancy]{Directional discrepancy in two dimensions.}

\author[D. Bilyk]{Dmitriy Bilyk }
\address {Department of Mathematics, University of South Carolina, Columbia, SC, 29208}
\email{bilyk@math.sc.edu}

\author[X. Ma]{Xiaomin Ma}
\address {Mathematics Department, Brown University, Providence, RI, 02912}
\email{xiaomin@math.brown.edu}

\author[J. Pipher]{Jill Pipher}
\address {Mathematics Department, Brown University, Providence, RI, 02912}
\email{jpipher@math.brown.edu}

\author[C. Spencer]{Craig Spencer}
\address {Department of Mathematics, Kansas State University, Manhattan, KS, 66506}
\email{cvs@math.ksu.edu}



\thanks{The authors acknowledge the support of the National Science Foundation (the first author --  NSF grant DMS-0801036, the second and the third authors -- NSF grant DMS-0901139). In addition, the first and the fourth authors were supported by the NSF grant 0635607 at the Institute for Advanced Study.}


\maketitle

\begin{abstract}
In the present paper, we study the {\em{geometric discrepancy}} with respect to families of rotated rectangles. The well-known extremal cases are the axis-parallel rectangles (logarithmic discrepancy) and rectangles rotated in all possible directions (polynomial discrepancy). We study several intermediate situations: lacunary sequences of directions, lacunary sets of finite order, and sets with small Minkowski dimension. In each of these cases, extensions of a lemma due to Davenport allow us to construct appropriate rotations of the integer lattice which yield small discrepancy.
\end{abstract}

\section{Introduction}

In the present paper we address the following two-dimensional question in the theory of irregularities of distribution. Let $\Omega \subset [0, \pi/2]$ be a set of directions. We consider the collection of rectangles pointing in the directions of $\Omega$:
\begin{equation}\label{e.aomega}
{\mathcal A}_\Omega = \{ \textup{rectangles $R$ :\, one side of $R$ makes angle $\phi \in \Omega$ with the $x$-axis} \}.
\end{equation}
Taking a set of $N$ points in the unit square, ${\mathcal P}_N \subset [0,1]^2$, we measure its discrepancy with respect to ${\mathcal A}_\Omega$:
\begin{equation}\label{e.domega}
D_\Omega ({\mathcal P}_N) = \sup_{R\in {\mathcal A}_\Omega, \, R\subset [0,1]^2} |D_\Omega ({\mathcal P}_N,R)|=\sup_{R\in {\mathcal A}_\Omega, \, R\subset [0,1]^2} \bigg| \# {\mathcal P}_N \cap R  -  N\cdot |R| \bigg|.
\end{equation}
We are interested in the behavior of the quantity 
\begin{equation}\label{e.dn}
D_\Omega(N) = \inf_{{\mathcal P}_N \subset [0,1]^2} D_\Omega({\mathcal P}_N).
\end{equation}
as $N$ goes to infinity, depending on the properties of $\Omega$. It is also of interest to consider suitable (e.g., $L^2$) averages in place of the supremum in (\ref{e.domega}). 

The motivation for this question comes from several classical results:
\begin{itemize}
\item In the case $\Omega = \{ 0\}$, i.e. ${\mathcal A}_\Omega$ is the set of axis-parallel rectangles we have
\begin{equation}\label{e.d0}
D_\Omega \approx \log N.
\end{equation}
Here, and throughout the paper, we use the notation $A\lesssim B$ meaning that there exists an absolute constant $C$, independent of $N$, such that $ A\leq CB$, and write $A\approx B$ if $A\lesssim B \lesssim A$. The lower bound in the estimate above is a celebrated theorem of W. Schmidt \cite{Schmidt}, while the upper bound goes back to a century-old result due to Lerch \cite{Lerch}. The inequalities above continue to hold when $\Omega$ is finite (This result is essentially contained in \cite{ChTr}).

\item When $\Omega = [0, \pi/2]$, i.e. ${\mathcal A}_\Omega$ consists of rectangles rotated in all possible directions, we have
\begin{equation}\label{e.dall}
N^{\frac14} \lesssim D_{\Omega} (N) \lesssim N^\frac14 \log^\frac12 N.
\end{equation}
Here both inequalities are due to J. Beck (\cite{Beck1}, \cite{Beck2}). 
\end{itemize}

We see that the behavior of $D_\Omega (N)$ in these two extreme situations differs  drastically. We would like to know what happens in the intermediate cases, how the geometry of $\Omega$ effects the discrepancy, and where is the threshold between the logarithmic and polynomial estimates.

In this work we look at particular examples: $\Omega$ being 1) a lacunary sequence of directions; 2) a lacunary set of finite order (for the definition of such sets and a brief discussion of their role in analysis see \S \ref{s.lacm}); or 3) a set with small upper Minkowski dimension, and prove the following theorem:
\begin{theorem}\label{t.main}$ $

1) Let $\Omega$ be a lacunary sequence. Then we have
\begin{equation}\label{e.dlac}
D_\Omega (N) \lesssim \log^3 N.
\end{equation}

2)  Let $\Omega$ be a lacunary set of order $M >1$. Then we have
\begin{equation}\label{e.dlacm}
D_\Omega (N) \lesssim \log^{2M+1} N.
\end{equation}

3) Assume $\Omega$ has upper Minkowski dimension $0\le d<1$. In this case,
\begin{equation}\label{e.dmink}
D_\Omega(N) \lesssim N^{\frac{\tau}{2(\tau + 1)}+ \varepsilon},
\end{equation}
for any $\varepsilon >0$, where $\tau = {\frac{2}{(1-d)^2} -2  }$.
\end{theorem}

We should point out that, in view of (\ref{e.dall}), the last part yields a new non-trivial estimate only if $d$ is small enough.

In addition, we complement this theorem with the following $L^2$-averaging estimates. Denote  ${\mathcal A}'_\Omega = \{ R\in {\mathcal A}_\Omega \,:\, R\subset [0,1]^2 \}$, or, alternatively, one may define ${\mathcal A}'_\Omega = \{ R\in {\mathcal A}_\Omega \,:\, \textup{diam}(R)\le 1 \}$ with $[0,1]^2$  viewed as a torus. We have

\begin{theorem}\label{t.mainl2} Let $\mu $ be any probability measure on ${\mathcal A}'_\Omega$. Then

1) If $\Omega$ is a lacunary sequence, there exists ${\mathcal P}\subset [0,1]^2$, $\#{ \mathcal P} =N$ such that
\begin{equation}\label{e.dlacl2}
\left( \int_{{\mathcal A}'_\Omega}  | D_\Omega ( {\mathcal P}, R )|^2 d\mu(R) \right)^\frac12 \lesssim \log^{\frac52} N.
\end{equation}

2)  If $\Omega$ is a lacunary set of order $M >1$, there exists ${\mathcal P}\subset [0,1]^2$, $\#{ \mathcal P}=N$ such that
\begin{equation}\label{e.dlacml2}
\left( \int_{{\mathcal A}'_\Omega}  | D_\Omega ( {\mathcal P}, R )|^2 d\mu(R) \right)^\frac12 \lesssim \log^{2M+\frac12} N.
\end{equation}

3) If $\Omega$ has upper Minkowski dimension $0\le d<1$,  there exists ${\mathcal P}\subset [0,1]^2$, $\# {\mathcal P}=N$ such that
\begin{equation}\label{e.dminkl2}
\left( \int_{{\mathcal A}'_\Omega}   | D_\Omega ( {\mathcal P}, R )|^2 d\mu(R)  \right)^\frac12 \lesssim N^{\frac{\tau}{2(\tau + 1)}+ \varepsilon},
\end{equation}
for any $\varepsilon >0$, where $\tau = {\frac{2}{(1-d)^2} -2  }$ satisfies $\tau <1$.

\end{theorem}

Comparing the first two parts of the above theorem to those of Theorem \ref{t.main}, we see a manifestation of the well-known discrepancy theory principle  that the $L^\infty$ (extremal) and $L^2$ (average) discrepancies differ by a factor of $\sqrt{\log N}$.  This effect can be best seen if one compares  (\ref{e.d0}) to the famous Roth's $L^2$ lower bound \cite{Roth1} of the order $\log^{1/2} N$ (which is sharp, \cite{Dav2}). In addition, the lower bound in (\ref{e.dall}) is known to be sharp in the $L^2$ sense \cite{BeCh2}.

In addition, we also address a `sibling' problem: studying the discrepancy with respect to collections ${\mathcal B}_{\Omega, k}$ of convex polygons in $[0,1]^2$ with at most $k$ sides whose normals point in the directions defined by $\Omega$ (cf. \cite{BeCh}, \cite{ChTr} for earlier results) and prove inequalities analogous to Theorems \ref{t.main} and \ref{t.mainl2} (see Theorems \ref{t.conv} and \ref{t.convl2} in the text).

The paper is organized as follows. The core of the paper is \S 2 -- here we obtain new diophantine inequalities which enable us to construct well-distributed sets. Section 3 describes how such inequalities can be translated into  upper  discrepancy estimates for one-dimensional sequences. In \S 4, we deduce our main Theorem \ref{t.main}, and \S 5 deals with bounds for  the $L^2$ discrepancy in these settings. In the text, $\log n$ stands for $\max \{ 1, \log_2 n\}$.


\section{Cassels-Davenport diophantine approximation arguments}

In the case $\Omega=\{0\}$, one of the standard ways of constructing an example of a point-set satisfying the upper bound of (\ref{e.d0}) involves rotating the lattice $N^{-\frac12} {\mathbb Z}^2$ by an angle $\alpha$ so that the slope $\tan \alpha$ is 
a {\em badly approximable} number, that is, for all $p \in \mathbb Z$, all $q \in \mathbb N$ we have 
\begin{equation}\label{e.ba}
\left| \tan \alpha  - \frac{p}{q} \right| \gtrsim \frac{1}{q^2}.
\end{equation}

When $\Omega$ is an arbitrary finite set, the construction relies on the following result of Davenport \cite{Dav} (which we state here in a particular case, relevant to our problem)
\begin{lemma}\label{l.dav}
Let $\Omega = \{\theta_1, \theta_2 ,..., \theta_k \} \subset [0,\pi/2]$. Then there exists $\alpha \in [0,\pi/2]$ so that $$\tan (\alpha - \theta_1) , ... , \tan (\alpha - \theta_k) $$ are all badly approximable.
\end{lemma}
This allows us to find a rotation, which has a badly approximable slope with respect to all chosen directions $\theta_j$. Davenport has, in fact, proven this fact for more general functions in place of  the tangent. However, the argument is essentially due to Cassels \cite{Cass} who proved a similar result earlier with $\tan (\alpha - \theta_k)$ replaced by $\alpha - \theta_k$.

Thus, analogs of the lemma above for infinite sets $\Omega$ may provide us with examples of low-discrepancy point distributions with respect to rotated rectangles. However, claiming ``badly approximable" in the conclusion is, perhaps,  too optimistic. Instead, we shall  obtain results, in which inequalities similar to (\ref{e.ba}) have the right-hand side somewhat smaller than $1/q^2$. This, in turn, will lead to larger discrepancy bounds. 

\subsection{General approach}\label{s.gen}

We first outline a general approach to the proof of statements akin to Lemma \ref{l.dav} extending the ideas of Cassels and Davenport.  Assume that for a certain choice of parameters $R(n)$, $| I_n |$, $c(n)$, depending on the set $\Omega$, a  proposition of the following type holds:

\begin{proposition}\label{p.gen}
Let $\Omega \subset [0, \pi/2]$. There exists a sequence of nested intervals $I_{0} \supset I_{1}\supset ... \supset I_n \supset ...$ in $[0,\pi/2]$ with 
$|I_n| \rightarrow 0$ such that for all 
$\alpha \in I_n$ and all $p,q \in \mathbb Z$ with  $R(n) \le q< R(n+1)$ we have, for all $\theta \in \Omega$: 
\begin{equation}\label{e.eq0}
\left| \tan \left(\alpha - \theta \right) -\frac{p}{q} \right| > \frac{c(n)}{q^2}.
\end{equation}
\end{proposition} 

This would of course imply that:

\begin{lemma}\label{l.gen}
There exist $\alpha\in  [0,\pi/2]$ and $C>0$ such that for all $\theta \in \Omega$, all $p\in \mathbb Z$, $q\in \mathbb N$ we have
\begin{equation}\label{e.gen}
\left| \tan \left(\alpha -\theta \right) -\frac{p}{q} \right| > \frac{C}{q^2 \, f( q)},
\end{equation}
where the function $f(q)$ is determined by the relation between $c(n)$ and $R(n)$.
\end{lemma}

To prove (\ref{e.eq0}), one proceeds inductively. At the  $n^{th}$ step, the set $\Omega$ is covered by at most $N_n$ intervals of length  $\delta_n$:  the dependence  between $N_n$ and $\delta_n$ is governed by the geometry of the set $\Omega$:
\begin{itemize}
\item $N=\textup{const}$, if $\Omega$ is finite;
\item $N \lesssim \log \frac1\delta$, if $\Omega$ is lacunary;
\item $N \lesssim \log^M \frac1\delta$, if $\Omega$ is lacunary of order $M$;
\item $N \le C_\varepsilon \left(\frac1\delta \right)^{d+\varepsilon}$, if $\Omega$ has upper Minkowski dimension $d$.
\end{itemize}

 Next, one has to choose  parameters $R(n)$, $| I_n |$, $c(n)$, $\delta_n$, $N_n$ so that they satisfy two inequalities, for an appropriately chosen constant $C$ (We initially restrict our range of $\alpha$ to, say, $[\alpha_0, \pi/2 - \alpha_0]$, so that, for all $\theta \in \Omega$, $\alpha- \theta \in [-\pi/2 +\alpha_0, \pi/2 - \alpha_0]$, where the derivative of tangent is bounded above by some $C>0$):

\begin{equation}\label{e.ineq1}
\frac{2 c(n)}{R^2(n)} + C( |I_{n-1}| +  \delta_n) < \frac{1}{R^2(n+1)}\,\,\,\,\,\,\,\,\textup{and}
\end{equation}

\begin{equation}\label{e.ineq2}
|I_{n-1}| - N_n \left( \frac{2 c(n)}{R^2(n)} + \delta_n \right) \ge (N_n +1) |I_n|.
\end{equation}

Indeed, assuming that $I_{n-1}$ is constructed, fix one of the chosen intervals $\Omega_n^k$ of length $\delta_n$. Suppose that the inequality (\ref{e.eq0}) doesn't hold for two sets of numbers $\alpha', \alpha'' \in I_n$, $\theta', \theta'' \in \Omega_n^k$,  $p',p'' \in \mathbb Z$, $R(n) \le q', q'' < R(n+1)$, then by (\ref{e.ineq1})
\begin{eqnarray*}
\left| \frac{p'}{q'} - \frac{p''}{q''} \right| & \le & \left| \frac{p'}{q'} - \tan \left( \alpha' - \theta' \right)\right| + \left| \frac{p''}{q''} - \tan \left( \alpha'' - \theta'' \right)\right| + \left| \tan \left( \alpha' - \theta' \right)  -  \tan \left( \alpha'' - \theta'' \right)\right|  \\
& \le & \frac{2 c(n)}{R^2(n) } + C (| \alpha'-\alpha'' | + | \theta' -\theta''|)  \le \frac{2 c(n)}{R^2(n)} + C ( |I_{n-1}| +  \delta_n) < \frac{1}{R^2(n+1)},
\end{eqnarray*}
which shows that $p'/q'=p''/q''$ (for otherwise they would have to differ by at least $\frac1{R^2(n+1)}$), i.e., there is at most one fraction $p_k/q_k$ with $R(n) \le q', q'' < R(n+1)$ for each $\Omega_n^k$ for which (\ref{e.eq0}) is violated. 

This implies that the inequality is true for $\alpha$ away from $$S_n=\bigcup_{k=1}^{N_n} \left\{ {\tan^{-1} \left\{ \left[   \frac{p_k}{q_k}- \frac{ c(n)}{R^2(n)} , \frac{p_k}{q_k} + \frac{ c(n)}{R^2(n)}\right] \right\}} + \Omega_n^k\right\}.$$
Obviously, $|S_n| \le N_n \left( \frac{2 c(n)}{R^2(n)} + \delta_n \right)$ and $I_{n-1} \setminus S_n$ consists of at most $N_n+1$ intervals. Thus, the validity of (\ref{e.ineq2}) proves that $I_{n-1}\setminus S_n$ contains at least one interval of length $|I_n|$.

In particular, for a finite set $\Omega$, to prove Davenport's lemma (Lemma \ref{l.dav}), one can choose the parameters $R(n)=R^n$, $c(n)=c$ (for some $R,c>0$), $\delta_n=0$, $N_n = \# \Omega$. The task of proving similar lemmata for sets $\Omega$ of different types is therefore reduced to the proper choice of these parameters. The details are taken up in subsequent subsections.

\subsection{Lacunary sequences}

We recall that a sequence $\Omega = \{\omega_n\}_{n=1}^\infty$  is called lacunary if $\omega_{n+1}/ \omega_n < A$ for some $A<1$. For simplicity, we shall consider the set $\Omega = \{ 2^{-k} \}_{k=1}^\infty$, however the argument easily extends to more general lacunary sequences. The main geometrical feature of this set for our purposes is the fact that it can be covered by $\log_2 (1/\delta)$ intervals of length $\delta$. We prove

\begin{lemma}\label{l.lac}
There exist $\alpha\in  [0,\pi/2]$ and $C>0$ such that for all $k\in \mathbb N$, all $p\in \mathbb Z$, $q\in \mathbb N$ we have
\begin{equation}\label{e.lac}
\left| \tan \left(\alpha - 2^{-k}\right) -\frac{p}{q} \right| > \frac{C}{q^2 \log^2 q}.
\end{equation}
\end{lemma}
The result of the lemma will follow from the following proposition similar to Proposition \ref{p.gen}:

\begin{proposition}\label{p.lac}
There exists a sequence of nested intervals  $I_{n_0} \supset I_{n_0 +1}\supset ... \supset I_n \supset ...$  with 
$$|I_n|= \delta (n+2)^{-(n+2)} \big( \log (n+2) \big)^{-(n+2)},$$ such that for all 
$\alpha \in I_n$ and all $p,q \in \mathbb Z$ with  
$n^{\frac{n}{2}} \big( \log n \big)^{\frac{n}{2}}\le q < (n+1)^{\frac{n+1}{2}} \big( \log (n+1) \big)^{\frac{n+1}{2}}$ 
we have
\begin{equation}\label{e.eq1}
\left| \tan \left(\alpha - 2^{-k}\right) -\frac{p}{q} \right| > \frac{c(n)}{q^2},
\end{equation}
where $c(n) = \frac{c}{(n+1)^2 \log^2 (n+1)} $ (for some  absolute constants $c$, $\delta>0$ and $n_0 \in \mathbb N$.)
\end{proposition}

Indeed, the proposition  implies that there exists $\alpha$ such that  for all $k\in \mathbb N$ we have 
\begin{equation}\label{e.lac1}
\left| \tan \left(\alpha - 2^{-k}\right) -\frac{p}{q} \right| > \frac{c'}{q^2 \log^2 q}
\end{equation}
 for $q \ge q_0 = (n_0)^{\frac{n_0+1}{2}} \big( \log (n_0+1) \big)^{\frac{n_0+1}{2}}$ and for some $c'>0$.

Now consider $q\le q_0$. Choose integer $r$, $1\le r \le q_0$ so that $qr\ge q_0$. Then, if $q \ge 2$,
\begin{eqnarray*}
\left| \tan \left(\alpha - 2^{-k}\right) -\frac{p}{q} \right| & = \left| \tan \left(\alpha - 2^{-k}\right) -\frac{pr}{qr} \right| >\frac{c'}{(qr)^2 \log^2 (qr)}\\ 
& > \frac{c'}{q_0^2 (1+\log q_0)^2} \frac{1}{q^2 \log^2 q} = \frac{c''}{q^2 \log^2 q}
\end{eqnarray*}
for some constant $c''>0$. The case $q=1$  (without the $\log$) is easy.

{\em Proof of Proposition \ref{p.lac}.} We restrict the range of $\alpha$ to $[0, \pi/3]$ so that $\alpha - 2^{-k} \in [-1, \pi/3] \subset  [-\pi/3, \pi/3]$, so that the derivatives of $\tan \left( \alpha -2^{-k} \right)$ satisfy $$1 \le \frac1{\cos^2 \left( \alpha -2^{-k} \right)}\le 4 .$$

We arbitrarily choose an initial interval $I_{n_0-1} \subset  [-\pi/3, \pi/3]$ with length  $$|I_{n_0-1}| = \varepsilon (n_0+1)^{-(n_0+1)} \big( \log (n_0+1) \big)^{-(n_0+1)},$$ where $\varepsilon$ is a small constant,  and proceed to construct the sequence inductively.

At the $n^{th}$ step we cover $\Omega$ by at most  $N_n = 2 (n+1) \log (n+1)$ intervals of length $\delta_n = 2^{-N_n} = (n+1)^{-2(n+1)}$. We now show that with this choice of parameters ($c(n) = \frac{c}{(n+1)^2 \log^2 (n+1)} $, $R(n) = n^{\frac{n}{2}} \big( \log n \big)^{\frac{n}{2}}$,  $|I_n|= \varepsilon (n+2)^{-(n+2)} \big( \log (n+2) \big)^{-(n+2)}$) the inequalities (\ref{e.ineq1}) and (\ref{e.ineq2}) hold for $n$ large enough.

Indeed, one easily verifies (\ref{e.ineq1}):
\begin{equation}
  \frac{2 c(n)}{n^{n} (\log n)^n}  + 4 ( |I_{n-1}| + \delta_n ) <  \frac1{(n+1)^{n+1}(\log (n+1))^{n+1}} =\frac{1}{R^2(n+1)},
\end{equation} 
for $c$, $\varepsilon$  small. Inequality (\ref{e.ineq2}) is slightly more subtle, as in this case both sides have roughly the same order of magnitude in $n$, so a little extra care should be given to constants. It is easy to see that, if $c\ll \varepsilon$ and $n$ is large, the left-hand side satisfies 
\begin{equation}
|I_{n-1}| - N_n \left( \frac{2 c(n)}{R^2(n)} + \delta_n \right) > 0.99 | I_{n-1}|,
\end{equation}
(we have $N_n \frac{2 c(n)}{R^2(n)} \approx | I_{n-1}|$ and  $N_n \delta_n \ll | I_{n-1} |$ for $n$ large)

\noindent On the other hand, for the right-hand side
\begin{eqnarray*}
(N_n+1) |I_n| & \le & \varepsilon  ( 2 (n+1) \log (n+1) +1 )  \times (n+2)^{-(n+2)} \big( \log (n+2) \big)^{-(n+2)}  \\
& \le & \varepsilon \cdot 2.5  \cdot (n+2)^{-(n+1)} \big( \log (n+2) \big)^{-(n+1)} \quad\textup{for $n$ large}\\
& \le & \varepsilon \cdot 2.5 \cdot (n+1)^{-(n+1)} \big( \log (n+1) \big)^{-(n+1)}  \left( 1+\frac{1}{n+1} \right)^{-(n+1)}\\
& \le & \varepsilon \cdot \frac{2.5}{2.7} \cdot  (n+1)^{-(n+1)} \big( \log (n+1) \big)^{-(n+1)} \quad\textup{for $n$ large}\\
& <  & 0.99 \cdot \varepsilon \cdot (n+1)^{-(n+1)}(\log (n+1))^{-(n+1)}  =0.99 |I_{n-1}|,
\end{eqnarray*}
where the second inequality from the bottom holds because $e>2.7$. Thus, (\ref{e.ineq2}) holds and the proof is finished.

\subsection{Lacunary sets of finite order}\label{s.lacm}
We now turn our attention to lacunary sets of finite order. They are defined inductively
\begin{definition}
Lacunary set of order one is a lacunary sequence. We call a set $\Omega$ lacunary of order $M$ if it can be
covered by the union of a lacunary set $\Omega'$ of order $M-1$ with lacunary sequences
converging to every point of $\Omega'$.
\end{definition}

These sets play an important role in analysis. In particular, recently M. Bateman \cite{Bat} proved that the directional maximal function
\begin{equation}
{\mathcal M}_\Omega f  (x) = \sup_{R \in {\mathcal A}_\Omega:\,\,x \in R} \frac{1}{|R|} \int_R |f(x)|\, dx,
\end{equation}
where  ${\mathcal A}_\Omega$ is as defined in (\ref{e.aomega}), is bounded on $L^p({\mathbb R}^2)$, $1<p<\infty$, if and only if $\Omega$ is covered by a finite union of lacunary sets of finite order. This condition is also equivalent to the fact that $\Omega$ does not ``admit Kakeya sets" (for details see \cite{Bat}, \cite{SS}).

One can check that a lacunary set of order $M$ can be covered by ${\mathcal{O}} (\log^M (1/\delta))$ intervals of length $\delta$. A simple example of a lacunary set of order $M$ is a set 
\begin{equation}\label{e.omlacm}
\Omega= \{ 2^{-j_1}+2^{-j_2}+...+2^{-j_M} \}_{j_1,...,j_M \in{\mathbb N}}.
\end{equation}
 In our setting, we have the following statement about such sets:

\begin{lemma}\label{l.lacm}
Let $\Omega \subset [0, \pi/2]$ be a lacunary set of order $M\ge 1$. Then there exist $\alpha\in  [0,\pi/2]$ and $C>0$ such that for all $\theta \in \Omega$, all $p\in \mathbb Z$, $q\in \mathbb N$ we have
\begin{equation}\label{e.lacm}
\left| \tan \left(\alpha - \theta \right) -\frac{p}{q} \right| > \frac{C}{q^2 \log^{2M} q}.
\end{equation}
\end{lemma}

This lemma is a generalization of Lemma \ref{l.lac}. For simplicity we deal with $\Omega$ as in (\ref{e.omlacm}) in which case $N(\delta) = \log_2^M (M/\delta)$. We follow the general approach of \S \ref{s.gen} and verify that inequalities (\ref{e.ineq1}) and (\ref{e.ineq2}) hold for the following choice of parameters 
\begin{eqnarray*}
R(n) & =  & (Mn)^{\frac{Mn}{2}} \big( \log n \big)^{\frac{Mn}{2}},\\
 |I_n| & = & \varepsilon (M(n+2))^{-M(n+2)} \big( \log (n+2) \big)^{-M(n+2)},\\ 
 c(n) & = & \frac{c}{(M(n+1))^{2M} \log^{2M} (n+1)},\\ 
 N_n & =  & (2M)^M (n+1)^M \log^M (n+1),\\ 
 \delta_n & = & M 2^{-N_n^{1/M}} = M(n+1)^{-2M(n+1)}.
 \end{eqnarray*}
 The proof is verbatim the same as that of Proposition \ref{p.lac}.

\subsection{Sets of fractional Minkowski dimension}
We now turn to  an analogous lemma for the case when the set of directions has non-negative upper Minkowski dimension. Recall that the upper Minkowski dimension of a set $\Omega \subset \mathbb R$ is defined as the infimum of exponents $d$ such that for any $0<\delta \ll 1$ the set $E$ can be covered by ${\mathcal O}(\delta^{-d})$ intervals of length $\delta$.

\begin{lemma}\label{l.mink}
Let $\Omega \subset (0, \pi/2)$ be a set of upper Minkowski dimension $d<1$. Then, for each $\varepsilon >0$, there exists $\alpha \in \mathbb R$ and a constant $c>0$ such that  for all $\gamma \in \Omega$, all $p\in \mathbb Z$, $q \in \mathbb Z_+$ we have
\begin{equation}\label{e.mink}
\left| \tan (\alpha-\gamma) - \frac{p}{q} \right| > c \,{q^{-\frac{2}{(1-d)^2}-\varepsilon}}.
\end{equation}
\end{lemma} 

The proof is again based on the approach described in \S \ref{s.gen}. Fix $t\in(d,1)$ and denote $a=\frac{1}{1-t}$. We shall construct a system of nested intervals $I_n$ with length $|I_n| = \varepsilon_1 2^{-2a^{n+2}}$ such that for $p\in \mathbb Z$, $R(n) = 2^{a^{n}} \le q < 2^{a^{n+1}}=R(n+1)$ we have, for all $\alpha \in I_{n}$, $$\left| \tan (\alpha-\theta) - \frac{p}{q} \right| > \frac{c(n)}{q^2}, $$ where $c(n)= c 2^{-2a^n (a^2 -1)}$. The lemma follows from this construction, since $c(n)\gtrsim q^{-2(a^2-1)}$ for this range of $q$'s.

Initially, restrict the attention to $\alpha$ in $(\alpha_0, \pi/2-\alpha_0)$, $\alpha_0 >0$, so that $\alpha -\theta$ stays away from $\pm \pi/2$ and the derivative of $\tan(\alpha-\theta)$ is bounded above by some $C>0$ in absolute value.

Assume $I_{n-1}$ is constructed and consider $2^{a^{n}} \le q < 2^{a^{n+1}}$. Now fix a number $s$ so that $d<s<t$. We cover $\Omega$ by  at most $N_n=C_s \delta_n^{-s}$ intervals of length $\delta_n = \varepsilon_2 2^{-2a^{n+2}}$. Inequality (\ref{e.ineq1}) is obviously satisfied
\begin{equation}
\frac{2 c(n)}{2^{2a^n}} + C ( |I_{n-1}| + \delta ) <  2^{-2a^{n+1}} = \frac1{R^2(n+1)},
\end{equation}
if the constants $c, \varepsilon_1, \varepsilon_2$ are small enough. 
\begin{eqnarray*}
N_n \cdot \left( \frac{ 2c(n)}{2^{2a^n}} +\delta_n\right) & \le &  C_s \delta_n^{-s} \left( \frac{ 2c(n)}{2^{2a^n}} +\delta_n \right)\\
& \le &  C_s \delta_n^{-t} \left( \frac{ 2c(n)}{2^{2a^n}} +\delta_n \right)\\
& = & C_s \left( 2 c\delta_n^{-t} 2^{-2a^{n+2}}  +  \delta_n^{1-t}   \right)\\
& = & C_s  \left( 2c \varepsilon_2^{-t} 2^{2a^{n+2} t}2^{-2a^{n+2}}  +  \varepsilon_2^{1-t} 2^{-2a^{n+1}}  \right)\\
& = & C_s  2^{-2a^{n+1}}  \left( 2c \varepsilon_2^{-t}   +  \varepsilon_2^{1-t}  \right)\\
& <  & \frac12 \varepsilon_1 2^{-2a^{n+1}} =\frac12  |I_{n-1}|,
\end{eqnarray*}
if $\varepsilon_2$ and $c$ are small (notice that $a(1-t)=1$).  Then $ | I_{n-1} | - N_n \cdot \left( \frac{ 2c(n)}{2^{2a^n}} +\delta_n\right)  \ge \frac12 | I_{n-1} |$ and
\begin{equation}
(C_s \delta^{-s} + 1)|I_{n}|  \lesssim  2^{2a^{n+2} s} 2^{-2a^{n+2}} = 2^{-2a^{n+2} (1-s)} = 2^{-2a^{n+1} \left(\frac{1-s}{1-t} \right)} \approx |I_{n-1}|^{\frac{1-s}{1-t}}.
\end{equation}
Since $\frac{1-s}{1-t}>1$, we conclude that $(C_s \delta^{-s} + 1)|I_{n}| < \frac12 |I_{n-1}| $ for $n$ large enough. Thus (\ref{e.ineq2}) holds and the proof is finished.

\section{One-dimensional discrepancy estimates}

Denote by $\| \theta \|$ the distance from $\theta$ to the nearest integer. We say that a real number $\theta$ is of type $<\psi$ for some non-decreasing function $\psi$ on $\mathbb R_+$ if for all natural $q$ we have $q \| q\theta \|  > 1/\psi(q)$, in other words for all $p\in \mathbb Z$, $q \in \mathbb N$ we have 
\begin{equation}\label{e.psi}
\left|  \theta - \frac{p}{q} \right|>\frac{1}{q^2\cdot \psi(q)}.
\end{equation}
In particular, our results in the previous section imply that the numbers $\tan(\alpha-\gamma) $ are of type $<\psi$ with
\begin{itemize}
\item $\psi(q) = C\, \log^2 q$ in the lacunary case,
\item $\psi(q) = C\, \log^{2M} q$ in the ``lacunary of order $M$" case,
\item $\psi(q) =C\, q^{\frac{2}{(1-d)^2} -2 + \varepsilon }$ in the case of upper Minkowski dimension $d$.
\end{itemize}

For a sequence ${\mathbf \omega} = \{ \omega_n \}_{n=1}^\infty \subset [0,1]$ its discrepancy is defined as 
\begin{equation}\label{e.d1d}
D_N ({\mathbf \omega} ) = \sup_{x\in [0,1]}  \bigg|  \# \big\{ \{\omega_1, ..., \omega_N \} \cap [0,x) \big\} - Nx \bigg|
\end{equation}
The Erd\"{o}s-Turan inequality (in a simplified form) says that, for any sequence $\omega \subset [0,1]$ 
\begin{equation}\label{e.e-t}
D_N(\omega)  \lesssim \frac{N}{m} + \sum_{h=1}^m \frac{1}{h} \left| \sum_{n=1}^N e^{2\pi i h \omega_n} \right|
\end{equation} for all natural numbers $m$. It is particularly convenient to apply it to the sequence of the form $\{n\theta\}$, since in this case 
$$\left| \sum_{n=1}^N e^{2\pi i h n \theta} \right|\le \frac{2}{|e^{2\pi i h \theta}-1|}= \frac{1}{|\sin(\pi h \theta)|} 
= \frac{1}{\sin(\pi \| h \theta \|)} \le  \frac{1}{2 \| h \theta \|},$$
since $\sin(\pi x) \ge 2x$ for $x\in [0, 1/2]$.
Thus, we obtain
\begin{equation}\label{e.e-t2}
D_N(\{n\theta\})  \lesssim \frac{N}{m} + \sum_{h=1}^m \frac{1}{h \| h \theta \|}.
\end{equation}

If the number $\theta$ is of type $<\psi$, then the last sum above can be estimated as follows
 (see e.g., Exercise 3.12, page 131, \cite{KN})
\begin{equation}\label{e.h2}
\sum_{h=1}^m \frac{1}{h \| h \theta \|} \lesssim  \log^2 m + \psi(m) + \sum_{h=1}^m \frac{\psi(h)}{h}.
\end{equation}

{\em Remark.} The proof of the estimate above is  somewhat delicate; a more straightforward summation by parts argument (Lemma 3.3, page 123, \cite{KN}) would have given
\begin{equation}\label{e.h1}
\sum_{h=1}^m \frac{1}{h \| h \theta \|} \lesssim \psi(2m) \log m + \sum_{h=1}^m \frac{\psi(2h) \log h}{h}.
\end{equation}
However, in the case of lacunary directions, this inequality would have given us a weaker bound. It is interesting to note that in the case $\psi=const$, i.e. $\theta$ is badly approximable, both estimates, (\ref{e.h2}) and (\ref{e.h1}), only yield $\log^2 N$ as opposed to the sharp $\log^1N$.

\begin{itemize}
\item The case $\psi(q) = C\, \log^2 q$. We have
 $$ \sum_{h=1}^m \frac{1}{h\| h \theta \|} \lesssim \log^2 m + \sum_{h=1}^m \frac{\log^2 h}{h} 
 \approx \log^3 m, $$ while (\ref{e.h1}) would only have given $\log^4 m$.
 Thus, for the discrepancy, inequality (\ref{e.e-t2}) with $m\approx N$ yields
 \begin{equation}\label{e.d1lac}
D_N(\{n\theta\}) \lesssim \log^3 N.
 \end{equation}
 
 \item More generally, in the case $\psi(q) = C\, \log^{2M} q$, we obtain
  \begin{equation}\label{e.d1lacm}
D_N(\{n\theta\}) \lesssim \log^{2M+1} N.
 \end{equation}

\item The case $\psi(q) =C\, q^{\frac{2}{(1-d)^2} -2 + \varepsilon }$. Denote $\tau = {\frac{2}{(1-d)^2} -2 + \varepsilon }$. From (\ref{e.h2}) we get $$ \sum_{h=1}^m \frac{1}{h \| h \theta \|} \lesssim  m^\tau + \sum_{h=1}^m h^{\tau-1} \approx m^\tau .$$ Inequality (\ref{e.e-t2}) with $m\approx N^{\frac{1}{\tau+1}}$ shows that the discrepancy satisfies
\begin{equation}\label{e.d1mink}
D_N(\{n\theta\}) \lesssim N^{\frac{\tau}{\tau+1}}.
\end{equation}
\end{itemize}

\section{Discrepancy with respect to rotated rectangles}

In the present section we demonstrate how one can translate the one-dimensional discrepancy estimates into the estimates for $D_\Omega (N)$. These ideas are classical and go  back to Roth \cite{Roth1}. The exposition of this and the next sections  essentially follows the papers of Beck and Chen \cite{BeCh} and Chen and Travaglini \cite{ChTr}.

The examples providing the upper bounds will be obtained using a rotation of the lattice $(N^{-1/2} {\mathbb Z})^2$. However, for technical reasons, it will be easier to rotate the unit square and the rectangles instead and leave the lattice intact. In addition, we shall consider a rescaled version of the problem.

Assume $\Omega$ is as described in parts 1,2, or 3 of Theorem \ref{t.main}. Let $\alpha$ be the angle  provided by Lemma \ref{l.lac}, \ref{l.lacm}, or \ref{l.mink}, respectively. Denote by $V$ the square $[0, N^{1/2})$ rotated counterclockwise by $\alpha$, and by ${\mathcal A}_{\Omega, \alpha}$ the family of all rectangles $R\subset  V$ which have a side that is either parallel to a side of $V$ or  makes angle $\theta - \alpha$ with the $x$-axis for some $\theta \in\Omega$ . (Strictly speaking, we should have applied Lemma \ref{l.lac},  \ref{l.lacm},  or \ref{l.mink} to the set $\Omega \cup \{0\} \cup (\Omega + \pi/2) \cup \{ \pi/2 \}$. It is easy to see that this change does not alter the proof.) For $R\subset V$, consider the quantity $D(R) = \# \{ {\mathbb Z}^2 \cap R \} - |R|$. We have the following lemma:
\begin{lemma}\label{l.main}$ $

1) Let $\Omega$ be a lacunary sequence. For any $R\in{\mathcal A}_{\Omega, \alpha}$  we have
\begin{equation}\label{e.dlac1}
D(R) \lesssim \log^3 N.
\end{equation}

2) Let $\Omega$ be a lacunary set of order $M$. For any $R\in{\mathcal A}_{\Omega, \alpha}$  we have
\begin{equation}\label{e.dlacm1}
D(R) \lesssim \log^{2M+1} N.
\end{equation}

3) Assume $\Omega$ has upper Minkowski dimension $0<d<1$. In this case, for each $R \in{\mathcal A}_{\Omega, \alpha}$,
\begin{equation}\label{e.dmink1}
D_\Omega(N) \lesssim N^{ \frac{\tau}{2 (\tau + 1)}+ \varepsilon},
\end{equation}
for any $\varepsilon >0$, where $\tau = {\frac{2}{(1-d)^2} -2  }$.
\end{lemma}
 
 We first show that the lemma above implies our main theorem.
 
 {\em Proof of Theorem \ref{t.main}.}   Denote by ${\mathcal P}_\alpha$ the intersection of the lattice $(N^{-1/2}{\mathbb Z} ) ^2$ rotated by $\alpha$ and $[0,1]^2$. The only obstacle to proving the theorem is the fact that  ${\mathcal P}_\alpha$ does not necessarily contain precisely $N$ points.  Let ${\mathcal P}'_\alpha$ be a set of $N$ points obtained from  ${\mathcal P}_\alpha$ by arbitrarily adding or removing $| \# {\mathcal P}_\alpha  - N | $ points. Let $F(N)$ stand for the right-hand side of the inequality we are proving ((\ref{e.dlac}), (\ref{e.dlacm}), or (\ref{e.dmink})). ``Unscaling" the estimates of Lemma \ref{l.main} and taking $R=[0,1]^2$, we obtain $$ | \# {\mathcal P}_\alpha  - N | \lesssim F(N). $$
 
 Then, for any $R\in {\mathcal A}_\Omega$ we have, again using Lemma \ref{l.main}
 \begin{eqnarray*}
 \bigg| \# {\mathcal P}'_\alpha \cap R  - N |R| \bigg| 
 & \le  &\bigg| \# {\mathcal P}_\alpha \cap R  - N |R| \bigg|  +  \bigg| \# {\mathcal P}_\alpha \cap R  -  \# {\mathcal P}'_\alpha \cap R  \bigg| \\
 & \lesssim & F (N)  +  | \# {\mathcal P}_\alpha  - N | \lesssim F(N),
 \end{eqnarray*}
 which finishes the proof. $\Box$
 
 {\em Remark.}  In view, of inequality (\ref{e.dall}), for any $\Omega$ we have the bound  $D_\Omega (N) \lesssim N^{1/4} \log^{1/2} N$. Thus, the bound arising from (\ref{e.dmink}) is meaningful only if $\frac12-\frac{1}{2(1+\tau)} < \frac14$, i.e. $\tau < 1$.  So, in the context of rotated rectangles, this estimate is interesting only if the set of rotations has low Minkowski dimension: 
 \begin{equation}\label{e.drestrict}
 d< 1-\left( \frac23 \right)^\frac12 \approx 0.1835... .
 \end{equation}

We now prove Lemma \ref{l.main}.  For each point ${\bf n} = (n_1,n_2) \in {\mathbb Z}^2$, consider a square of area one centered around it $$S({\bf n}) = \left[ n_1 - \frac12, n_1+\frac12\right) \times \left[ n_2 - \frac12, n_2+\frac12\right).$$
Obviously, we can write: $$ D(R) = \sum_{{\bf n} \in {\mathbb Z}^2} D(R\cap S({\bf n})). $$ Denote the sides of $R$ by $T_1$, $T_2$, $T_3$, $T_4$. Set $${\mathcal N}^- =\{ {\bf n}:\, S({\bf n}) \cap T_i =  \emptyset, \,\,\textup{for all } i=1,2,3,4\}, $$ i.e. the set of those $\bf n$ for which the corresponding square lies entirely within or entirely outside $R$ -- for such squares $ D(R\cap S({\bf n}))=0$.

Also, take  $${\mathcal N}^+  =\{ {\bf n}:\, S({\bf n}) \cap T_i \neq \emptyset, S({\bf n}) \cap T_{i+1} \neq \emptyset, \,\,\textup{for some } i=1,2,3,4\}, $$ (the addition is mod $4$) to be those $\bf n$ for which $S({\bf n})$ contains a corner of $R$. We have $\# {\mathcal N}^+ \le 4$ and $| D(R\cap S({\bf n})) | \le 1$, thus $\sum_{{\bf n} \in {\mathcal N}^+} D(R\cap S({\bf n})) \le 4$.

Finally, for $i=1,2,3,4$, set $$ {\mathcal N}^i = \{ {\bf n}:\, S({\bf n}) \cap T_i \neq \emptyset, \,\textup{but}\, {\bf n} \not\in {\mathcal N}^+ \}, $$ to be the centers of  those squares which intersect  the side $T_i$ but do not contain any corners.
The collections ${\mathcal N}^i$ are not necessarily disjoint, e.g., when $R$ is a thin rectangle. However, we have the following useful fact:
\begin{proposition}\label{p.square}
Let $R$ be a convex polygon with sides $T_1$,...,$T_m$. Denote by $T_j^*$ the halfplane with boundary $T_j$ which contains $R$.  Assume the square $S({\bf n})$ intersects $R$ but does not contain any vertices of $R$.  Let $T_{j_1}$, ... , $T_{j_k}$ be the sides of $R$ that intersect $S({\bf n})$. Then
\begin{equation}\label{e.square}
D( R \cap S({\bf n})) = \sum_{i=1}^k D(T_{j_i}^* \cap S({\bf n})).
\end{equation}
\end{proposition}

We use the fact that discrepancy is an additive measure and that $D(S({\bf n}))=0$. Then
$$0 = D(S({\bf n})) = \sum_{i=1}^k \bigg( D(S({\bf n})) -  D(T_{j_i}^* \cap S({\bf n})) \bigg) + D( R \cap S({\bf n})). \,\,\,\square
$$

 Since ${\mathbb Z}^2 =  {\mathcal N}^- \cup {\mathcal N}^+ \cup {\mathcal N}^1 \cup...\cup {\mathcal N}^4$, it remains to estimate the terms  $\sum_{{\bf n} \in {\mathcal N}^j} D(T_j^* \cap S({\bf n}))$. Assume that the $j^{th}$ side of $R$ lies on the line $\tan \phi = \frac{y_2-a_2}{y_1 -a_1}$, i.e. 
 $$y_2 = a_2 +  (y_1 - a_1) \tan \phi$$
  for some constants $a_1$, $a_2$ and $\phi = \alpha - \theta$ or  $\phi =  \alpha - \theta +\pi/2$. Let $I_j= \{ n_1 \in {\mathbf Z}: \, (n_1, n_2)\in {\mathcal N}^j \, \textup{for some }\, n_2\in{\mathbb Z}\}$ be the projection of the ${\mathcal N}^j$ onto the $x$-axis. Fix $n\in I_j$ and let $h\in \mathbb Z$ be the smallest number such that $(n,h)\in {\mathcal N}^j$. Then it is easy to see that (here we assume that $R$ is below $T_j$, the other case is analogous)
\begin{equation}\label{e.points}
\sum_{{\bf n} \in {\mathcal N}^j,\, n_2=n} \# \{ {\mathbb Z}^2 \cap T_j^* \cap S({\bf n})\} = [ y_2 (n_1) - h  +1 ],
\end{equation}
and the area of the trapezoid is 
\begin{equation}\label{e.area}
\sum_{{\bf n} \in {\mathcal N}^j,\, n_2=n} | T_j^*\cap S({\bf n}) | =  y_2 (n_1) - h  + \frac12 .
\end{equation}
(This relation may fail when $n$ is an endpoint of $I_j$, but this gives us a bounded error.) Thus, the discrepancy can be described by the ``sawtooth" function, $\psi(x) =  x - [x] - \frac12 = \{x\} - \frac12$,
\begin{equation}\label{e.side}
\sum_{{\bf n} \in {\mathcal N}^j} D(T_j^*\cap S({\bf n}))  = \pm   \sum_{n\in I_j}  \psi ( c - n \tan \phi ) .
\end{equation}

The ``sawtooth" function arises naturally in one dimensional discrepancy. If we define, for a sequence $\omega$, 
$$D_N(\omega, x) = \bigg|  \# \big\{ \{\omega_1, ..., \omega_N \} \cap [0,x) \big\} - Nx \bigg| , $$
one can easily check that  
\begin{equation}\label{e.dpsi}
D_N(\omega, x) = \sum_{n=1}^N  \bigg( \psi (\omega_n - x) - \psi (\omega_n) \bigg). 
\end{equation}

Since $x \in [0,1]$ is arbitrary, it is possible to show that for all $x\in [0,1]$
\begin{equation}\label{e.psid}
 \left| \sum_{n=1}   \psi (\omega_n - x) \right| \lesssim D_N(\omega).
 \end{equation}
Indeed,  one can find a point $x\in [0,1]$ with $ D_N(\omega, x) = \sum_{n=1}^N   \psi (\omega_n) $ (see, e.g., the proof of Erd\"{o}s-Turan in \cite{KN}), thus $\left| \sum_{n=1}^N   \psi (\omega_n) \right| \le D_N(\omega)$, but then for any $x\in[0,1]$, 
$\left| \sum_{n=1}  \psi (\omega_n - x) \right| \le 2 D_N(\omega) $. Thus, (\ref{e.side}) and (\ref{e.psid}) imply
\begin{equation}\label{e.dd}
\left| \sum_{{\bf n} \in {\mathcal N}^j} D(R\cap S({\bf n})) \right| \lesssim D_{| I_j |} (\omega).
\end{equation}

Obviously, $|I_j| \lesssim N^\frac12$. This fact, together with inequality (\ref{e.dd}) and the results of the previous section, proves the lemma. $\Box$

To conclude this section, we formulate analogous results on the discrepancy with respect to convex polygons. We omit the proofs as they are verbatim the same as the proof of the main theorem.

Let $\Omega$ be a set of directions. Denote by ${\mathcal B}_{\Omega, k}$ the collection of all convex polygons in $[0,1]^2$ with at most $k$ sides whose normals belong to $\pm \Omega$ and set  $$D_{\Omega,k} (N) = \inf_{{ \mathcal P}_N :\, \# { \mathcal P}_N = N}  \sup_{B\in {\mathcal B}_{\Omega, k}} \bigg| \# {\mathcal P}_N \cap B -  N\cdot |B| \bigg|.$$
The following theorem holds (notice that the implied constants depend on $k$):
\begin{theorem}\label{t.conv}$ $

1) Let $\Omega$ be a finite union of lacunary sets of order at most $M\ge1$. Then we have
\begin{equation}\label{e.cdlac}
D_{\Omega, k} (N) \lesssim_k \log^{2M+1} N.
\end{equation}

2) Assume $\Omega$ has upper Minkowski dimension $0<d<1$. In this case,
\begin{equation}\label{e.cdmink}
D_{\Omega,k} (N) \lesssim_k N^{\frac{\tau}{2(\tau + 1)}+ \varepsilon},
\end{equation}

for any $\varepsilon >0$, where $\tau = {\frac{2}{(1-d)^2} -2  }$.
\end{theorem}

\section {An upper bound for the  $L^2$ discrepancy}

We now prove Theorem \ref{t.mainl2}. In this case, the point set with low $L^2$ discrepancy is given by a suitably  {\emph {shifted}} rotation of the lattice $(N^{-1/2} {\mathbb Z})^2$; the idea of using random shifts to obtain distributions with low average discrepancy was first introduced by Roth \cite{Roth2}.    As in the previous section we consider a rescaled and rotated version of the problem, that is  we set $V$ to be the square$[0,{N^{1/2}}]^2$ rotated counterclockwise by the angle $\alpha$ given by the  Lemma \ref{l.lac}, \ref{l.lacm}, or \ref{l.mink}.  Assume $ {\mathcal A}_{\Omega, \alpha}$ is the family of all rectangles $R\subset  V$ which have a side that is either parallel to a side of $V$ or  makes angle $\theta - \alpha$ with the $x$-axis for some $\theta \in\Omega$ and fix a rectangle $R\in  {\mathcal A}_{\Omega, \alpha} $.

 For any $\omega \in [0,1 ]^2$ define the shift of the integer lattice ${\mathbb Z}_{\omega}^2=\omega+{\mathbb Z}^2$. Consider the quantity $D_{\omega} (R)=D({\mathbb Z}_{\omega}^2, R)=\# \{{\mathbb Z}_{\omega}^2\cap R\}- |R|$. We estimate  the mean square of the shifted discrepancies  in the following lemma:

\begin{lemma}\label{l.fourier}$ $

1) Let $\Omega$ be a lacunary set of order $M\ge 1$. For any $R\in  {\mathcal A}_{\Omega, \alpha}$, we have
\begin{equation}\label{e.fourier1}
\int_{[0,1]^2} \mid D({\mathbb Z}_{\omega}^2, R) \mid ^2 d{\omega} \lesssim \log^{4M+1} N
\end{equation}

2) Let $\Omega$ be a set of upper Minkowski dimension $d<1$. For any $R\in  {\mathcal A}_{\Omega, \alpha}$, we have
\begin{equation}\label{e.fourier2}
\int_{[0,1]^2} \mid D({\mathbb Z}_{\omega}^2, R) \mid^2 d{\omega} \lesssim N^{ \frac{\tau}{ \tau + 1}+ \varepsilon},
\end{equation}
for any $\varepsilon >0$, where $\tau = {\frac{2}{(1-d)^2} -2  }$ and $\tau < 1$.
\end{lemma}

The lemma relies on the following important calculation which goes back to Davenport \cite{Dav2}  (see also Beck and Chen \cite{BeCh}).  Recall that $\parallel x \parallel=\min_{n \in Z}
\mid x-n \mid$ denotes the distance from $x$ to the nearest integer. We have
\begin{lemma}\label{l.Fouriersum} Let $I$ be a finite  interval of consecutive integers.

1) Assume $\tan \phi$ satisfies $\nu \| \nu \tan \phi \| > \frac{c } {\log^{2M} \nu}$, for all $\nu \in \mathbb N$ . Then
\begin{equation}
\sum _{\nu = 1}^{\infty}{\frac{1} {\nu^2}} \left| {{\sum_{n \in {I}}} e^{-2\pi i \nu n \tan \phi}  } \right|^2 \lesssim \log^{4M+1} {\vert{I}\vert}.
\end{equation}

2) Assume $\tan \phi$ satisfies $\nu \| \nu \tan \phi \| > c \nu^{-\tau+\varepsilon}$, for all $\varepsilon>0$, where $0\le \tau <1$. Then
\begin{equation}
\sum _{\nu = 1}^{\infty}{\frac{1} {\nu^2}} \left| {{\sum_{n \in {I}}} e^{-2\pi i \nu n \tan \phi}  } \right|^2 \lesssim {\vert{I}\vert}^{\frac{2\tau}{\tau+1}+\varepsilon'}, \,\,\,\textup{where}\,\,\varepsilon' ={\mathcal O}(\varepsilon).
\end{equation}
\end{lemma}

{\emph {Proof.}} We will use a simple  fact that
$$ \left| \sum_{n \in I}  e^{-2\pi i \nu n \tan \phi}  \right| \, \lesssim \min \{\vert {I}\vert, \parallel \nu  \tan \phi \parallel^{-1}\}.$$
We deal with part one first:
\begin{eqnarray*}
S & = &\sum _{\nu = 1}^{\infty}{\frac{1} {\nu^2}} \left| {{\sum_{n \in {I}}} e^{-2\pi i \nu n \tan \phi}  } \right|^2 \\
& \lesssim & \sum_{h=1}^{\infty}{2^{-2h}}\sum _{2^{h-1} \le \nu <2^h} \min\{\vert {I}\vert^2, \parallel  \nu \tan \phi \parallel^{-2}\}.
\end{eqnarray*}
Notice that our assumption on $\tan \phi $ implies that if $ 2^{h-1} \le \nu  <2^h$, then $\parallel  \nu \tan \phi \parallel> {c \over {2^h h^{2M}}}$. On the other hand, for any pair $h, p \in {\mathbb N}$, there are at most two values of $\nu$ satisfying $2^{h-1} \le \nu  <2^h$ and $p{c \over {2^h h^{2M}}} \le \parallel \nu \tan \phi \parallel < (p+1){c \over {2^h h^{2M}}}$. Indeed,  otherwise the difference $(\nu_1-\nu_2)$ of two of them would contradict
the assumption. We have

\begin {eqnarray*}
S & \lesssim &  {\sum_{h=1}^ \infty}{\sum_{p=1}^ \infty}\min\{ 2^{-2h} \vert {I}\vert^2, p^{-2}h^{4M}\}\\
& = & {\sum_{2^h \le \mid{I} \mid}} {\sum_{p=1}^\infty} \min\{ 2^{-2h} \vert {I}\vert^2, p^{-2}h^{4M}\}
+  {\sum_{2^h > \mid{I} \mid}} {\sum_{p=1}^\infty}\min\{ 2^{-2h} \vert {I}\vert^2, p^{-2}h^{4M}\} \\
& \lesssim & {\sum_{2^h \le \mid{I} \mid}} {\sum_{p=1}^\infty} \,p^{-2}h^{4M} + {\sum_{2^h > \mid{I} \mid}} \left(2^{-2h}\vert {I} \vert ^2 2^h \vert {I} \vert ^{-1}h^{2M} + {\sum_{p > 2^h h^{2M} \mid {I} \mid^{-1}}}h^{4M}p^{-2} \right) \\
& \lesssim &{\sum_{2^h \le \mid{I} \mid}}h^{4M}+{\sum_{2^h > \mid{I} \mid}}2^{-h}\vert{I} \vert h^{2M} \\
& \lesssim & \log^{4M+1} \mid{I}\mid.
\end{eqnarray*}

Part 2 is proved in a similar fashion. The choice of  $\phi$ yields that, for $ 2^{h-1} \le \nu <2^h$, we have $\parallel \nu \tan \phi \parallel> {c2^{h(-1-\tau-\varepsilon)}}$. And as before, for any pair $h, p \in {\mathbb N}$, no more than two values of $\nu$ satisfy $2^{h-1} \le \nu <2^h$ and $p{c2^{h(-1-\tau-\varepsilon)}} \le \parallel \nu \tan \phi \parallel < (p+1){c2^{h(-1-\tau-\varepsilon)}}$. Thus

\begin {eqnarray*}
S & \lesssim & {\sum_{h=1}^ \infty}{\sum_{p=1}^ \infty} \min\{ 2^{-2h} \vert {I}\vert^2, p^{-2}2^{2h(\tau+\varepsilon)}\}\\
& = & {\sum_{2^{h(1+\tau)} \le \mid{I} \mid}} {\sum_{p=1}^\infty}\min\{ 2^{-2h} \vert {I}\vert^2, p^{-2}2^{2h(\tau+\varepsilon)}\}
+{\sum_{2^{h(1+\tau)} > \mid{I} \mid}} {\sum_{p=1}^\infty}\min\{ 2^{-2h} \vert {I}\vert^2, p^{-2}2^{2h(\tau+\varepsilon)}\} \\
& \lesssim & {\sum_{2^{h(1+\tau)} \le \mid{I} \mid}} {\sum_{p=1}^\infty}\, p^{-2}2^{2h(\tau+\varepsilon)}+{\sum_{2^{h(1+\tau)} > \mid{I} \mid}}\left(    2^{-2h}\vert {I} \vert ^2 {2^{h(1+\tau+\varepsilon)}\mid {I} \mid^{-1} + {\sum_{p>2^{h(1+\tau+\varepsilon)}\mid {I} \mid^{-1}}}p^{-2}2^{2h(\tau+\varepsilon)}}\right) \\
& \lesssim & {\sum_{2^{h(1+\tau)} \le \mid{I} \mid}}2^{2h(\tau+\varepsilon)}+{\sum_{2^{h(1+\tau)} > \mid{I} \mid}}2^{h(-1+\tau+\varepsilon)}\vert{I} \vert \\
& \lesssim & \mid I \mid ^{{2\tau \over {1+\tau}}+\varepsilon'},
\end{eqnarray*}
where $\tau<1$ is required for the second sum in the penultimate line above to converge for any choice of $\varepsilon>0$.

We turn to the proof of  Lemma \ref{l.fourier}. For any ${\bf n}=(n_1,n_2) \in {\mathbb Z}^2$, ${\bf \omega}=(\omega_1,\omega_2) \in [0,1]^2$, define
$$S({\bf n},{\bf \omega})=[n_1+\omega_1-1/2, n_1+\omega_1+1/2) \times [n_2+\omega_2-1/2, n_2+\omega_2+1/2).$$

\noindent Also define
${\mathcal N}^+=\{{\bf n}: \exists {\bf \omega'} \in [0,1]^2 \,\,\textup{such that } \,\, S({\bf n},{\bf \omega'}) \textup{ contains a vertex of $R$ }\},$
and
\begin{eqnarray*}
{\mathcal N}=\{ {\bf n}: &\exists {\bf \omega'} \in [0,1]^2 \,\,\textup{ such that }\, S({\bf n},{\bf \omega'}) \cap R \ne \emptyset,
\,\textup{ and} \\ 
& \forall {\bf \omega} \in [0,1]^2, \,\,\,\,S({\bf n},{\bf \omega}) \,\,\textup{ contains no vertex of $R$ }  \}.
\end{eqnarray*}

Let $\widetilde{\mathcal N}={\mathcal N}^+ \cup {\mathcal N}^-$. Then one can see that $D_\omega(R)= \sum_{{\bf n} \in \widetilde{\mathcal N}} D_\omega(R \cap S({\bf n},{\bf \omega}))$. Obviously, $\# {\mathcal N}^+ = {\mathcal O}(1)$ and it remains to deal with ${\mathcal N}$. Write ${\mathcal N}= {\mathcal N}^1 \cup...\cup {\mathcal N}^4$ in a natural way. Using Proposition \ref{p.square}, we can rewrite the discrepancy
\begin{equation}
\sum_{{\bf n} \in {\mathcal N}} D_\omega(R \cap S({\bf n},{\bf \omega}))=\sum_{j=1}^4 \sum_{{\bf n} \in {\mathcal N}^j} D_\omega(S({\bf n},{\bf \omega}) \cap T^*_j)
\end{equation}
where $T^*_{j}$ is the halfplane defined by the $j^{th}$ side of $R$ (see Proposition \ref{p.square}).

For each $j=1,...,4$, define $I_j=\{n_1 \in {\mathbf Z}:\,  \exists n_2 \,\textup{such that}\, (n_1,n_2) \in {\mathcal N}^j\}$. Applying the argument, similar to the one preceding (\ref{e.side}), we express the discrepancy arising from the $j^{th}$ side in terms of  the ``sawtooth" function $\psi(x)$, up to a bounded error:
\begin{equation}\label{e.sidew}
\sum_{{\bf n} \in {\mathcal N}^j} D_\omega (S({\bf n},{\bf \omega}) \cap T^*_j) = \pm \sum_{n_1\in {I_j}} \psi (a_2-\omega_2+(n_1-a_1+\omega_1) \tan \phi)
\end{equation}

The ``sawtooth" function $\psi(x)$ has the Fourier expansion
$-\sum_{\nu \ne 0} {{e^{(2 \pi i\nu x)}} \over {2 \pi i \nu}}$. Hence, using Parseval's theorem, one easily obtains
\begin{equation}
\int_{[0,1]^2} \left| \sum_{{\bf n} \in {\mathcal N}^j} D_\omega (S({\bf n},{\bf \omega}) \cap T^*_j) \right|^2 d\omega \lesssim \sum _{\nu = 1}^{\infty}{\frac{1} {\nu^2}} \left| {{\sum_{n \in {I_j}}} e^{-2\pi i \nu n \tan \phi}  } \right|^2.
 \end{equation}

By applying Lemma \ref{l.Fouriersum} and the fact that, for each $j$, we have $\mid I_j \mid ={\mathcal O}(N^{1/2})$, we finish the proof of Lemma
\ref{l.fourier}.

We are now ready to prove Theorem \ref{t.mainl2}. Let $\mu $ be any probability measure on ${\mathcal A}'_\Omega$ and consider the induced probability measure $\mu'$ on the set ${\mathcal A}_{\Omega,\alpha}$ of rectangles $R\subset V$ (see the beginning of this section). Since, by Lemma \ref{l.fourier},  $$ \int_{[0,1]^2} \mid D({\mathbb Z}_{\omega}^2, R) \mid ^2 d{\omega} \lesssim F(N),$$
(where $F(N)$ denotes the right-hand side of (\ref{e.fourier1})  or (\ref{e.fourier2}), respectively), it follows that there exist $\omega_0 \in [0,1]^2$ such that $$\int_{{\mathcal A}_{\Omega,\alpha}}  \mid D({\mathbb Z}_{\omega_0}^2, R) \mid ^2 \, d\mu' (R) \lesssim F(N).$$
The only obstacle to finishing the proof is the fact that ${\mathbb Z}_{\omega_0}^2 \cap V$ does not necessarily contain precisely $N$ points. However, this can be handled as explained in the proof of Theorem \ref{t.main}.

{\emph{ Remark 1.}} Part 2 of Lemma \ref{l.fourier} required that $\tau<1$, which yields the same restriction $d<1-(2/3)^{1/2} \approx 0.1835.....$ that arises in the $L^ \infty$ case, (\ref{e.drestrict}), for a different reason.

{\emph{ Remark 2.}} Often, when considering $L^2$ averages, it is more convenient, instead of imposing the condition $R\subset [0,1]^2$,  to deal with all rectangles $R\in {\mathcal A}_\Omega$ with $\textup{diam} (R) \le 1$, while treating $[0,1]^2$ as a torus. In this case, the proof of Theorem \ref{t.mainl2} presented above undergoes only minor changes: modulo $V$, any rectangle $R$ with $\textup{diam} (R) \le N^{1/2}$ can be represented as at most 4 polygons  contained in    $V$, having  at most 6 sides each.

{\emph{ Remark 3.}} It is easy to see that the same argument also applies to convex polygons
with a bounded number of sides. Thus we also have the following theorem.

Let, as before,  ${\mathcal B}_{\Omega, k}$ denote the collection of all convex polygons in $[0,1]^2$ with at most $k$ sides whose normals belong to $\pm \Omega$ and set, for ${\mathcal P}\subset [0,1]^2$ with $\# {\mathcal P} = N$ and for $B\in   {\mathcal B}_{\Omega, k}$, $$D_{\Omega,k} ({\mathcal P}, B) =  \bigg| \# {\mathcal P} \cap B -  N\cdot |B| \bigg|.$$

\begin{theorem}\label{t.convl2} Let $\sigma$ be any probability measure on ${\mathcal B}_{\Omega, k}$

1) Let $\Omega$ be a finite union of lacunary sets of order at most $M\ge1$. Then there exists ${\mathcal P}\subset [0,1]^2$ with $\# {\mathcal P} = N$ such that
\begin{equation}\label{e.cdlacl2}
\left( \int_{{\mathcal B}_{\Omega, k}} |  D_{\Omega, k} ({\mathcal P}, B) |^2 \, d\sigma(B) \right)^\frac12 \lesssim_k \log^{2M+\frac12} N.
\end{equation}

2) Assume $\Omega$ has upper Minkowski dimension $0\le d<1$. In this case, there exists ${\mathcal P}\subset [0,1]^2$ with $\# {\mathcal P} = N$ such that
\begin{equation}\label{e.cdminkl2}
\left( \int_{{\mathcal B}_{\Omega, k}} |  D_{\Omega, k} ({\mathcal P}, B) |^2 \, d\sigma(B) \right)^\frac12 \lesssim_k N^{\frac{\tau}{2(\tau + 1)}+ \varepsilon},
\end{equation}

for any $\varepsilon >0$, where $\tau = {\frac{2}{(1-d)^2} -2  }$ satisfies $\tau <1$.
\end{theorem}

{\em{ Acknowledgments.}} All four authors have been supported by the National Science Foundation. The first three authors are sincerely grateful to the American Institute of Mathematics for the warm welcome during the workshop ``Small ball inequalities in analysis, probability, and irregularities of distribution".  Dmitriy Bilyk and Craig Spencer would like to thank the Institute for Advanced Study  for hospitality. Craig Spencer was also supported by the NSA Young Investigators Grant. Jill Pipher and Dmitriy Bilyk would like to express their gratitude to Centre de Recerca Matem\`{a}tica for support. In addition, the authors are indebted to William Chen and Giancarlo Travaglini for numerous interesting  and fruitful discussions.





\end{document}